\documentclass[a4,12pt,epsfig]{article}
\usepackage{graphicx}
\usepackage{subfig}
\usepackage{amsmath}
\usepackage{amssymb}
\usepackage{amsfonts}
\usepackage{euscript}
\begin{document}
\bibliographystyle{plain}
\newcommand{\bea}{\begin{eqnarray}}
\newcommand{\eea}{\end{eqnarray}}
\newcommand{\bfmN}{{\mbox{\boldmath{$N$}}}}
\newcommand{\bfmx}{{\mbox{\boldmath{$x$}}}}
\newcommand{\bfmv}{{\mbox{\boldmath{$v$}}}}
\newcommand{\se}{\setcounter{equation}{0}}
\newtheorem{corollary}{Corollary}[section]
\newtheorem{example}{Example}[section]
\newtheorem{definition}{Definition}[section]
\newtheorem{theorem}{Theorem}[section]
\newtheorem{proposition}{Proposition}[section]
\newtheorem{lemma}{Lemma}[section]
\newtheorem{remark}{Remark}[section]
\newtheorem{result}{Result}[section]
\newcommand{\vtwo}{\vskip 4ex}
\newcommand{\vthree}{\vskip 6ex}
\newcommand{\vfour}{\vspace*{8ex}}
\newcommand{\hone}{\mbox{\hspace{1em}}}
\newcommand{\hon}{\mbox{\hspace{1em}}}
\newcommand{\htwo}{\mbox{\hspace{2em}}}
\newcommand{\hthree}{\mbox{\hspace{3em}}}
\newcommand{\hfour}{\mbox{\hspace{4em}}}
\newcommand{\von}{\vskip 1ex}
\newcommand{\vone}{\vskip 2ex}
\newcommand{\n}{\mathfrak{n} }
\newcommand{\m}{\mathfrak{m} }
\newcommand{\q}{\mathfrak{q} }
\newcommand{\aF}{\mathfrak{a} }

\newcommand{\kl}{\mathcal{K}}
\newcommand{\p}{\mathcal{P}}
\newcommand{\Lt}{\mathcal{L}}
\newcommand{\bv}{{\mbox{\boldmath{$v$}}}}
\newcommand{\bc}{{\mbox{\boldmath{$c$}}}}
\newcommand{\bx}{{\mbox{\boldmath{$x$}}}}
\newcommand{\br}{{\mbox{\boldmath{$r$}}}}
\newcommand{\bs}{{\mbox{\boldmath{$s$}}}}
\newcommand{\bb}{{\mbox{\boldmath{$b$}}}}
\newcommand{\ba}{{\mbox{\boldmath{$a$}}}}
\newcommand{\bn}{{\mbox{\boldmath{$n$}}}}
\newcommand{\bp}{{\mbox{\boldmath{$p$}}}}
\newcommand{\by}{{\mbox{\boldmath{$y$}}}}
\newcommand{\bz}{{\mbox{\boldmath{$z$}}}}
\newcommand{\be}{{\mbox{\boldmath{$e$}}}}
\newcommand{\proof}{\noindent {\sc Proof :} \par }
\newcommand{\bP}{{\mbox{\boldmath{$P$}}}}

\newcommand{\M}{\mathcal{M}}
\newcommand{\R}{\mathbb{R}}
\newcommand{\Q}{\mathbb{Q}}
\newcommand{\Z}{\mathbb{Z}}
\newcommand{\N}{\mathbb{N}}
\newcommand{\C}{\mathbb{C}}
\newcommand{\xar}{\longrightarrow}
\newcommand{\ov}{\overline}
 \newcommand{\rt}{\rightarrow}
 \newcommand{\om}{\omega}
 \newcommand{\wh}{\widehat }
 \newcommand{\wt}{\widetilde }
 \newcommand{\g}{\Gamma}
 \newcommand{\lm}{\lambda}

\newcommand{\eN}{\EuScript{N}}
\newcommand{\ncom}{\newcommand}
\newcommand{\norm}{\|\;\;\|}
\newcommand{\inp}[2]{\langle{#1},\,{#2} \rangle}
\newcommand{\nrm}[1]{\parallel {#1} \parallel}
\newcommand{\nrms}[1]{\parallel {#1} \parallel^2}
\title{On the convergence of quasilinear viscous approximations with BV initial data}
\author{ Ramesh Mondal\footnote{ramesh@math.iitb.ac.in}\,\, and S. Sivaji Ganesh\footnote{siva@math.iitb.ac.in}}
\maketitle{}
\begin{abstract}
We show that the almost everywhere limit of quasilinear viscous approximations is the unique entropy
solution (in the sense of {\it Bardos-Leroux-Nedelec}) of the corresponding scalar conservation laws on a bounded domain 
in $\mathbb{R}^{d}$ whenever the initial data is essentially bounded and a function of bounded variation.
 \end{abstract}
\section{Introduction} 
Let $\Omega$ be a bounded domain in $\mathbb{R}^{d}$ with smooth boundary $\partial \Omega$. For $T >0$, denote $\Omega_{T}:= 
\Omega\times(0,T)$. We write the initial boundary value problem $\left(\mbox{IBVP}\right)$ for scalar conservation laws given by
 \begin{subequations}\label{ivp.cl}
\begin{eqnarray}
  u_t + \nabla\cdot f(u) =0& \mbox{in }\Omega_T,\label{ivp.cl.a}\\
u(x,t)= 0&\mbox{on}\,\,\partial \Omega\times(0,T),\label{ivp.cl.b}\\
  u(x,0) = u_0(x)& x\in \Omega,\label{ivp.cl.c}
  \end{eqnarray}
\end{subequations}
where $f=(f_{1},f_{2},\cdots,f_{d})$ is the flux function and $u_{0}$ is the initial condition.\\
Denote by $f_{\varepsilon}$, the regularizations of the flux function $f=(f_{1},f_{2},\cdots,f_{d})$ by the sequence of mollifiers 
$\rho_{\varepsilon}$ defined on $\R^d$. It is given by 
$f_{\varepsilon} := \left(f_{1\varepsilon},f_{2\varepsilon},\cdots,f_{d\varepsilon}\right)$, where
\begin{eqnarray*}\label{regularized.eqn2}
 f_{j\varepsilon}:= f_{j}\ast\tilde\rho_{\varepsilon} \, (j=1,2,\cdots,d)
\end{eqnarray*}
\vspace{0.2cm}
Consider the IBVP for regularized generalized viscosity problem 
\begin{subequations}\label{regularized.IBVP}
\begin{eqnarray}
 u^\varepsilon_{t} + \nabla \cdot f_{\varepsilon}(u^{\varepsilon}) = \varepsilon\,\nabla\cdot
 \left(B(u^\varepsilon)\,\nabla u^\varepsilon\right)&\mbox{in }\Omega_{T},\label{regularized.IBVP.a} \\
    u^\varepsilon(x,t)= 0&\,\,\,\,\mbox{on}\,\, \partial \Omega\times(0,T),\label{regularized.IBVP.b}\\
u^{\varepsilon}(x,0) = u_{0\varepsilon}(x)& x\in \Omega,\label{regularized.IBVP.c}
\end{eqnarray}
\end{subequations}
indexed by $\varepsilon>0$. Let us now give the hypothesis on $f,\,B,\,u_{0}$ and $u_{0\varepsilon}$. \\
\vspace{0.2cm}\\
\noindent{\bf Hypothesis E:}
\begin{enumerate}
 \item Let $f\in C^1(\R)$, $f^\prime\in L^\infty(\R)$, and denote 
 $$\|f^\prime\|_{L^\infty(\R)}:=\sup_{y\in\R}|f^\prime(y)|.$$
 \item Let $B\in C^3(\R)\cap L^\infty(\R)$, and there exists an $r>0$ such that $B\geq r$.
 \item  We denote all those elements of $L^{\infty}({\Omega})$ whose essential 
 support are compact subsets of $\Omega$ by $L^{\infty}_c({\Omega})$ and let $\mbox{BV}(\Omega)$ be the space of functions of bounded
 variations in $\Omega$. Let $u_{0}$ be in $ \mbox{BV}(\Omega)\cap L^{\infty}_{c}(\Omega)$ and denote by $u_{0\varepsilon}$, the 
 regularizations of the initial data $u_{0}$ by the sequence of mollifiers $\tilde{\rho}_{\varepsilon}$ 
 defined on $\R$, {\it i.e.}, $u_{0\varepsilon}:= u_{0}\ast\tilde{\rho}_{\varepsilon}$. Denote $I:=[-\|u_0\|_{\infty},\|u_0\|_{\infty}]$.
\end{enumerate}
\noindent{\bf Hypothesis F}
\begin{enumerate}
 \item Let $f\in \left(C^1(\R)\right)^d$, $f^\prime\in \left(L^\infty(\R)\right)^d$, and denote 
 $$\|f^\prime\|_{\left(L^\infty(\R)\right)^d}:=\max_{1\leq j\leq d}\,\sup_{y\in\R}|f^\prime_j(y)|.$$
 \item Let $B\in C^3(\R)\cap L^\infty(\R)$, and there exists an $r>0$ such that $B\geq r$.
 \item  We denote the set of all infinitely differentiable functions with compact 
 support in $\Omega$ by $\mathcal{D}(\Omega)$ and denote $W^{1,1}_{0}(\Omega):=\overline{\mathcal{D}(\Omega)}$ in $W^{1,1}(\Omega)$. 
 Let $u_{0}\in W^{1,1}_{0}(\Omega)\cap C(\overline{\Omega})$. Let $u_{0\varepsilon}$ be in $\mathcal{D}(\Omega)$ such that 
 for all $\varepsilon>0$, there exists a constant $A >0$ such that $\|u_{0\varepsilon}\|_{L^{\infty}(\Omega)}\leq A$ and 
 $u_{0\varepsilon}\to u_0$  in $W^{1,1}_{0}(\Omega)$ as $\varepsilon\to 0$. Denote $I:=\left[-A, A\right]$.
\end{enumerate}
The aim of this article is to prove that the {\it a.e.} limit of sequence of solutions 
$\left(u^{\varepsilon}\right)$ to \eqref{regularized.IBVP}(called quasilinear viscous approximations) is the unique entropy solution
for IBVP \eqref{ivp.cl}. In this context, we have two main results (Theorem \ref{paper3.BVestimates.theorem2} and 
Theorem \ref{paper3.BVestimates.theorem3}) depending on the regularity of the initial data. In \cite{Ramesh}, we have proved that 
the {\it a.e.} limit of quasilinear viscous approximations is the unique entropy solution in the sense of {\it Bardos et.al} whenever
the initial data lies in $W^{1,\infty}_{c}(\Omega)$ and in \cite{Mondal}, we showed that the {\it a.e.} limit of quasilinear viscous 
approximations is the unique entropy solution in the sense of {\it Bardos et.al} whenever the initial data lies in 
$H^{1}(\Omega)\cap L^{\infty}_{c}(\Omega)$ using compensated compactness. In this article, we are able to show that the {\it a.e.} limit
of quasilinear viscous approximations is the unique entropy solution in the sense of {\it Bardos et.al} whenever the initial data lies in 
$BV(\Omega)\cap L^{\infty}_{c}(\Omega)$ (see Hypothesis E) and we state our first result. 
\begin{theorem}\label{paper3.BVestimates.theorem2}
{\rm Let $f,\,B,\,u_{0}$ and $u_{0\varepsilon}$ satisfy Hypothesis E. Then the {\it a.e.} limit of the quasilinear viscous approximations 
 $\left(u^{\varepsilon}\right)$ satisfying \eqref{regularized.IBVP} is the unique entropy solution of IBVP \eqref{ivp.cl} in the sense 
 of {\it Bardos et.al} \cite{MR542510}.}
\end{theorem}
For the case artificial viscosity problem,{\it i.e.},$B(\cdot)\equiv 1$, it is enough to assume initial data in 
$BV(\Omega)\cap L^{\infty}_{c}(\Omega)$ for eshtablishing the BV estimates of quasilinear viscous approximations 
$\left(u^{\varepsilon}\right)$. But in the case of regularized viscosity problem, we are unable to achieve $L^{1}-$ estimate of the time 
derivative of $\left(u^{\varepsilon}\right)$ using usual technique \cite{MR542510}, \cite{MR1304494} because of the presence of 
non-constant $B$. In \cite{Ramesh}, we eshtablish BV estimate with initial data in $W^{1,\infty}_{c}(\Omega)$ and in 
\cite{Mondal}, we eshtablish $L^{1}-$ estimate of time derivative with initial data in $H^{1}(\Omega)\cap L^{\infty}_{c}(\Omega)$. 
But in the present article, we are able to give an alternate proof of the $L^{1}-$estimate of time derivative of quasilinear viscous 
approximations $\left(u^{\varepsilon}\right)$ with initial data in $BV(\Omega)\cap L^{\infty}_{c}(\Omega)$.\\
\vspace{0.2cm}\\
In \cite{Ramesh}, we removed compact essential support of initial data by showing the {\it a.e.} limit of quasilinear viscous approximations
is the unique entropy solution in the sense of {\it Bardos et.al} whenever the initial data lies in 
$H^{1}_{0}(\Omega)\cap C(\overline{\Omega})$. In this article we show that the {\it a.e.} limit of quasilinear viscous approximations is 
the unique entropy solution in the sense of {\it Bardos et.al} whenever the initial data lies in 
$W^{1,1}_{0}(\Omega)\cap C(\overline{\Omega})$ (see Hypothsis F). We now state our second result. 
\begin{theorem}\label{paper3.BVestimates.theorem3}
{\rm Let $f,\,B,\,u_{0}$ and $u_{0\varepsilon}$ satisfy Hypothesis F. Then the {\it a.e.} limit of the quasilinear viscous approximations 
 $\left(u^{\varepsilon}\right)$ satisfying \eqref{regularized.IBVP} is the unique entropy solution of IBVP \eqref{ivp.cl} in the sense 
 of {\it Bardos et.al} \cite{MR542510}.}
\end{theorem}
In this article, the main difficulty is to eshtablish the BV estimate of quasilinear viscous approximations with initial data as given in 
Hypothsis E and Hypothesis F. Then the proof of Theorem \ref{paper3.BVestimates.theorem2} and Theorem \ref{paper3.BVestimates.theorem3} follow 
from \cite{Ramesh}.\\
\vspace{0.1cm}\\
The plan of the paper is the following. In Section 2, we prove a few properties of quasilinear viscous approximations and in Section 3, we
eshtablish the BV estimate of quasilinear viscous approximations $\left(u^{\varepsilon}\right)$ and prove 
Theorem \ref{paper3.BVestimates.theorem2} and Theorem \ref{paper3.BVestimates.theorem3}.
\section{Properties of quasilinear viscous approximations}
We prove the following result for the existence and uniqueness, maximum principle of quasilinear viscous approximations 
$\left(u^{\varepsilon}\right)$.   
\begin{lemma}\label{regularized.BVestimates.theorem3}
Let $f, B, u_{0}$ and $u_{0\varepsilon}$ satisfy Hypothesis E. Then there exists a unique solution of \eqref{regularized.IBVP} in 
$C^{4+\beta,\frac{4+\beta}{2}}(\overline{\Omega_{T}})$, for every $0<\beta<1$ and the following estimate
  \begin{eqnarray}
  \|u^{\varepsilon}\|_{L^{\infty}(\Omega)}&\leq& \|u_{0}\|_{L^{\infty}(\Omega)}\,\,{\it a.e.}\,\,t\in(0,T)  \label{regularized.eqn19}
  \end{eqnarray}
 holds. 
\end{lemma}
In order to prove Lemma \ref{regularized.BVestimates.theorem3}, we use higher regularity and maximum principle of solutions to generalized 
viscosity probelm from \cite{Ramesh}. For generalized viscosity problem and for Hypotheses on $f,\,B,\, u_{0}$, we refer the reader to 
\cite[p.1]{Ramesh} and \cite[p.2]{Ramesh} respectively. We now state higher regularity result of generalized viscosity probelm  from 
\cite[p.18]{Ramesh}. This higher regularity result will be used to prove Lemma \ref{regularized.BVestimates.theorem3}.
\begin{theorem}[higher regularity]\label{chapHR85thm5}
Let $f,B,u_0$ satisfy Hypothesis A of \cite[p.2]{Ramesh}. Then the solutions of the IBVP for generalized viscosity problem belong to 
the space $C^{4+\beta,\frac{4 + \beta}{2}}(\overline{\Omega_{T}})$. Further,  $u_{tt}^{\varepsilon}\in C(\overline{\Omega_{T}})$.
\end{theorem}
We now state the following maximum principle for solutions of generalized viscosity problem from \cite[p.12]{Ramesh}.
\begin{theorem}[Maximum principle]\label{chap3thm1}
Let $f,\,B$ and $u_{0}$ satisfy Hypothesis A of \cite[p.2]{Ramesh}. Then any solution $u$ of generalized viscosity
problem in $C^{4+\beta,\frac{4 + \beta}{2}}(\overline{\Omega_{T}})$ satisfies the bound
\begin{equation}\label{eqnchap303}
||u^{\varepsilon}(\cdot,t)||_{L^{\infty}(\Omega)}\hspace{0.1cm}\leq\hspace{0.1cm}||u_{0}||_{L^{\infty}(\Omega)}\hspace{0.1cm}a.e.
\,\,t\in(0,T).
\end{equation}
\end{theorem}
\textbf{Proof of Lemma \ref{regularized.BVestimates.theorem3}:} Applying higher regularity result of solutions to generalized viscosity 
problem, {\it i.e.}, Theorem \ref{chapHR85thm5} to \eqref{regularized.IBVP}, we conclude the existence and uniquness
of quasilinear viscous approximations $\left(u^{\varepsilon}\right)$ to the regularized viscosity problem \eqref{regularized.IBVP}. \\
An application of Maximum principle, {\it i.e.}, Theorem \ref{chap3thm1} to quasilinear viscous approximations 
$\left(u^{\varepsilon}\right)$ and using $\left\|u_{0\varepsilon}\right\|_{L^{\infty}(\Omega)}\leq \|u_{0}\|_{L^{\infty}
(\Omega)}$, we obtain \eqref{regularized.eqn19}. This completes the proof of Lemma \ref{regularized.BVestimates.theorem3}.
\vspace{0.3cm}\\
Following exactly the same argument as in the proof of Lemma \ref{regularized.BVestimates.theorem3}, we conclude the following result.
\begin{lemma}\label{regularized.BVestimates.theorem3BB}
Let $f, B, u_{0}$ and $u_{0\varepsilon}$ satisfy Hypothesis F. Then there exists a unique solution of \eqref{regularized.IBVP} in 
$C^{4+\beta,\frac{4+\beta}{2}}(\overline{\Omega_{T}})$, for every $0<\beta<1$ and the following estimate
  \begin{eqnarray}
  \|u^{\varepsilon}\|_{L^{\infty}(\Omega)}&\leq& \|u_{0}\|_{L^{\infty}(\Omega)}\,\,{\it a.e.}\,\,t\in(0,T)  \label{regularized.eqn19BB}
  \end{eqnarray}
 holds. 
\end{lemma}
\vspace{0.2cm}
Applying Theorem 4.2 from \cite[p.30]{Ramesh} to quasilinear viscous approximations $\left(u^{\varepsilon}\right)$ as asserted in
Lemma \ref{regularized.BVestimates.theorem3}, we obtain
\begin{theorem}\label{Compactness.lemma.1}
Let $f,\,\,B,\,\,u_{0}$ and $u_{0\varepsilon}$ satisfy Hypothesis E. Let $u^{\varepsilon}$ be the unique solution to regularized 
viscosity problem \eqref{regularized.IBVP}. Then 
 \begin{eqnarray}\label{uniformnot.compactness.eqn1a}
 \displaystyle\sum_{j=1}^{d} \hspace{0.1cm}\left(\sqrt{\varepsilon}\Big\| \frac{\partial u^{\varepsilon}}{\partial x_{j}}\Big\|_{L^{2}
 (\Omega_{T})}\right)^{2} \leq\frac{1}{2r}\|u_{0\varepsilon}\|^{2}_{L^{2}(\Omega)}\leq\frac{1}{2r}
 \|u_{0}\|^{2}_{L^{\infty}(\Omega)}\,\,\mbox{Vol}(\Omega).
 \end{eqnarray}
\end{theorem}
\vspace{0.2cm}
Again applying Theorem 4.2 from \cite[p.30]{Ramesh} to quasilinear viscous approximations $\left(u^{\varepsilon}\right)$ as asserted in
Lemma \ref{regularized.BVestimates.theorem3BB}, we obtain
\begin{theorem}\label{Compactness.lemma.1D}
Let $f,\,\,B,\,\,u_{0}$ and $u_{0\varepsilon}$ satisfy Hypothesis F. Let $u^{\varepsilon}$ be the unique solution to regularized 
viscosity problem \eqref{regularized.IBVP}. Then 
 \begin{eqnarray}\label{uniformnot.compactness.eqn1aD}
 \displaystyle\sum_{j=1}^{d} \hspace{0.1cm}\left(\sqrt{\varepsilon}\Big\| \frac{\partial u^{\varepsilon}}{\partial x_{j}}\Big\|_{L^{2}
 (\Omega_{T})}\right)^{2} \leq\frac{1}{2r}\|u_{0\varepsilon}\|^{2}_{L^{2}(\Omega)}\leq\frac{1}{2r}\,\mbox{Vol}(\Omega)\, A^{2}.
 \end{eqnarray}
\end{theorem}
\section{BV estimates}
In this section we establish uniform $L^{1}(\Omega_{T})$ estimates of the first order derivatives of quasilinear viscous approximations 
$\left(u^{\varepsilon}\right)$ with respect to both time and space. We now state the BV estimates result.
\subsection{BV estimate with $u_{0}$ in $ BV(\Omega)\cap L^{\infty}_{c}(\Omega)$}
\begin{theorem}\label{paper3.BVestimate.Maintheorem}
 Let $f,\,B,\,u_{0}$ and $u_{0\varepsilon}$ be as in Hypothesis E and $\left(u^{\varepsilon}\right)$ be as asserted in 
 Lemma \ref{regularized.BVestimates.theorem3}. Then
 \begin{enumerate}
  \item for all $\varepsilon > 0$, the following inequality 
  \begin{eqnarray}\label{BVestimate.BVinitialdata.ean1}
   \|\nabla u^{\varepsilon}\|_{\left(L^{1}(\Omega_{T})\right)^{d}}\leq TV_{\Omega}(u_{0})
  \end{eqnarray}
  holds.
\item for all $\varepsilon >0$, the following inequality 
\begin{eqnarray}\label{BVestimate.BVinitialdata.ean2}
\left\|\frac{\partial u^{\varepsilon}}{\partial t}\right\|_{L^{1}(\Omega_{T})}&\leq& 2\|B^{\prime}\|_{L^{\infty}(I)}
 \frac{\mbox{Vol}(\Omega)}{2r}\|u_{0}\|_{L^{\infty}(\Omega)}^{2} + 2\displaystyle\max_{1\leq j\leq d}\left(\displaystyle\sup_{y\in I}
 \left|f_{j}^{\prime}(y)\right|\right)\,TV_{\Omega}(u_{0})\nonumber\\
 && + 2\,\|u_{0}\|_{L^{\infty}(\Omega)}\mbox{Vol}(\Omega) 
\end{eqnarray}
holds.
 \end{enumerate}
 Further there exists a subsequence $(u^{\varepsilon_{k}})$ of $(u^{\varepsilon})$, and a function $u\in L^{1}(\Omega_{T})$ such that 
\begin{eqnarray}\label{BVestimateeqn25}
u^{\varepsilon_{k}}\to u\,\,\mbox{in}\,L^1(\Omega_T),\label{BVestimateeqn25.a}\\
 u^{\varepsilon_{k}}\to u\,\,{\it a.e.}\,(x,t)\in \Omega_{T}\label{BVestimateeqn25.b}
\end{eqnarray}
as $k\to\infty$.
\end{theorem}
We now introduce signum function which will be used in the proof of Theorem \ref{paper3.BVestimate.Maintheorem}. For $n\in\mathbb{N}$, 
let $sg_{n}:\mathbb{R}\to\mathbb{R}$ be the sequence of functions given by 
$$
sg_{n}(s)=
       \begin{cases}
       1 & \,\mbox{if}\, s >\frac{1}{n},\\
       ns & \,\mbox{if}\, |s|\leq\frac{1}{n},\\
       -1 & \,\mbox{if}\, s < -\frac{1}{n},
       \end{cases}
$$
which converges pointwise to the signum function $sg:\mathbb{R}\to\mathbb{R}$ defined by
$$
sg(s)=
       \begin{cases}
       1 & \,\mbox{if}\, s > 0,\\
       0 & \,\mbox{if}\, s=0,\\
       -1 & \,\mbox{if}\, s < 0.
       \end{cases}
$$
The next result follows from \cite[p.67]{MR1304494} and is useful in proving Theorem \ref{paper3.BVestimate.Maintheorem}. 
\begin{lemma}\label{regularized.lem1}
 Let $u_{0}\in BV(\Omega)\cap L^{\infty}_{c}(\Omega)$ and $u_{0\varepsilon}$ be as in Hypothesis E. Then $u_{0\varepsilon}$ satisfies 
 the following bounds
  \begin{eqnarray}
    \|u_{0\varepsilon}\|_{L^{\infty}(\Omega)}\leq \|u_{0}\|_{L^{\infty}(\Omega)}\label{regularized.max.eqn1}\\[2mm]
     \|\nabla u_{0\varepsilon}\|_{\left(L^{1}(\Omega)\right)^{d}}\leq TV_{\Omega}(u_{0}) \label{regularized.max.eqn1a}
  \end{eqnarray}
There exists a constant $C> 0$ such that for all $\varepsilon > 0$, $u_{0\varepsilon}$ satisfies 
  \begin{eqnarray}\label{regularized.max.eqn1b}
   \|\Delta u_{0\varepsilon}\|_{L^{1}(\Omega)}\leq \frac{C}{\varepsilon}TV_{\Omega}(u_{0}).
  \end{eqnarray}
\end{lemma}
\textbf{Proof of Theorem \ref{paper3.BVestimate.Maintheorem}:} We prove Theorem \ref{paper3.BVestimate.Maintheorem} in three steps. In Step 1,
we show \eqref{BVestimate.BVinitialdata.ean1}, in Step 2, we show \eqref{BVestimate.BVinitialdata.ean2} and in Step 3, we show 
\eqref{BVestimateeqn25.a} and \eqref{BVestimateeqn25.b}.\\
\vspace{0.2cm}\\
\textbf{Step 1:} Applying Step 2  in the proof of BV estimate Theorem 4.1 from \cite[p.22]{Ramesh} with $f=f_{\varepsilon}$ and 
$u_{0}=u_{0\varepsilon}$, we arrive at
\begin{eqnarray}\label{BVestimate.BVinitialdata.ean3}\
 \|\nabla u^{\varepsilon}\|_{\left(L^{1}(\Omega)\right)^{d}}= \|\nabla u_{0\varepsilon}\|_{\left(L^{1}(\Omega)\right)^{d}}
\end{eqnarray}
Then applying Lemma \ref{regularized.lem1}, we conclude \eqref{BVestimate.BVinitialdata.ean1}.\\
\textbf{Step 2:} From equation \eqref{regularized.IBVP.a} of the regularized viscosity problem, we get
\begin{eqnarray}\label{BVestimate.BVinitialdata.ean4}
 u^\varepsilon_{t}  = \varepsilon\,B(u^\varepsilon)\,\Delta u^\varepsilon +\varepsilon\displaystyle\sum_{j=1}^{d}B^{\prime}(u^{\varepsilon})
 \left(\frac{\partial u^{\varepsilon}}{\partial x_{j}}\right)^{2} - \displaystyle\sum_{j=1}^{d}f_{j}^{\prime}(u^{\varepsilon})\frac{\partial 
 u^{\varepsilon}}{\partial x_{j}}\,&\mbox{on }\overline{\Omega_{T}}
\end{eqnarray}
From equation \eqref{BVestimate.BVinitialdata.ean4}, we get
\begin{eqnarray}\label{BVestimate.BVinitialdata.ean5}
 \varepsilon\,\int_{\Omega_{T}}B(u^\varepsilon)\,\left|\Delta u^\varepsilon\right|\,dx\,dt &\leq& \varepsilon\,\displaystyle\sum_{j=1}^{d}
 \int_{\Omega_{T}} B^{\prime}(u^{\varepsilon})\left(\frac{\partial u^{\varepsilon}}{\partial x_{j}}\right)^{2}\,dx\,dt + \displaystyle
 \sum_{j=1}^{d}\int_{\Omega}\left|f_{j}^{\prime}(u^{\varepsilon})\frac{\partial 
 u^{\varepsilon}}{\partial x_{j}}\right|\,dx\,dt\nonumber\\
 && + \int_{\Omega_{T}}\left|\frac{\partial u^\varepsilon}{\partial t}\right|\,dx\,dt\nonumber\\
 &\leq& \|B^{\prime}\|_{L^{\infty}(I)}\,\varepsilon\,\displaystyle\sum_{j=1}^{d}
 \int_{\Omega_{T}}\left(\frac{\partial u^{\varepsilon}}{\partial x_{j}}\right)^{2}\,dx\,dt\nonumber\\ &&+ \displaystyle\max_{1\leq j\leq d}
 \left(\displaystyle\sup_{y\in\mathbb{\mathbb{R}}}\left|f_{j}^{\prime}(y)\right|\right)
\displaystyle\sum_{j=1}^{d}\int_{\Omega}\left|\frac{\partial u^{\varepsilon}}{\partial x_{j}}\right|\,dx\,dt
 + \int_{\Omega_{T}}\left|\frac{\partial u^\varepsilon}{\partial t}\right|\,dx\,dt\nonumber\\
 {}
\end{eqnarray}
We now compute 
\begin{eqnarray}\label{BVestimate.BVinitialdata.ean6}
 \int_{\Omega_{T}}\left|\frac{\partial u^\varepsilon}{\partial t}\right|\,dx\,dt &=&\displaystyle\lim_{n\to\infty}\int_{\Omega_{T}}sg_{n}\left(\frac{\partial 
 u^\varepsilon}{\partial t}\right)\frac{\partial u^\varepsilon}{\partial t}\,dx\,dt
\end{eqnarray}
Using integration by parts, we get
\begin{eqnarray}\label{BVestimate.BVinitialdata.ean7}
 \int_{\Omega_{T}}\left|\frac{\partial u^\varepsilon}{\partial t}\right|\,dx\,dt &=&-\displaystyle\lim_{n\to\infty}\int_{\Omega_{T}}
 sg^{\prime}_{n}\left(\frac{\partial u^\varepsilon}{\partial t}\right)\,\frac{\partial^{2}u^{\varepsilon}}{\partial t^{2}}
  u^\varepsilon\,dx\,dt\nonumber\\
  &&+ \displaystyle\lim_{n\to\infty}\int_{\Omega}\left( sg_{n}\left(\frac{\partial }{\partial t}u^\varepsilon(x,T)\right)u^\varepsilon(x,T)
  - sg_{n}\left(\frac{\partial }{\partial t}u^\varepsilon(x,0)\right)u_{0\varepsilon}(x)\right)\,dx\nonumber\\
  {}
\end{eqnarray}
We want to show that
\begin{eqnarray}\label{BVestimate.BVinitialdata.ean8}
 \displaystyle\lim_{n\to\infty}\int_{\Omega_{T}}
 sg^{\prime}_{n}\left(\frac{\partial u^\varepsilon}{\partial t}\right)\,\frac{\partial^{2}u^{\varepsilon}}{\partial t^{2}}
  u^\varepsilon\,dx\,dt=0
\end{eqnarray}
The technique that we use to show \eqref{BVestimate.BVinitialdata.ean8}, was used by us in Step 2 in the proof of BV estimate 
Theorem 4.1 in \cite[p.25]{Ramesh}. We adapt here the argument.\\
 Denote 
$$A^{\varepsilon}:=\left\{(x,t)\in\Omega_{T}\,\,:\,\,\frac{\partial u^{\varepsilon}}{\partial t}=0\right\}.$$
Since $\frac{\partial u^{\varepsilon}}{\partial t}\in C^{1}(\overline{\Omega_{T}})$, $for \,\,{\it a.e.}\,\,(x,t)\in A^{\varepsilon}$,
using Stampacchia's theorem (see \cite{kesavan}), we conclude  
$\nabla_{x,t} \left(\frac{\partial u^{\varepsilon}}{\partial t}\right)=0$. In particular, 
$\frac{\partial^{2}u^{\varepsilon}}{\partial t^{2}}=0\,\,{\it a.e.}\,\,(x,t)\in A^{\varepsilon}$ and we have 
\eqref{BVestimate.BVinitialdata.ean8}. If $\Omega_{T}\smallsetminus A^{\varepsilon}=\varnothing$, then \eqref{BVestimate.BVinitialdata.ean8} 
follows trivially. Assume that $\Omega_{T}\smallsetminus A^{\varepsilon}\neq\varnothing$. For each 
$(x,t)\in\Omega_{T}\smallsetminus A^{\varepsilon}$, we have
\begin{eqnarray}\label{eqnarray}\label{BVestimate.BVinitialdata.ean9}
  sg^{\prime}_{n}\left(\frac{\partial u^\varepsilon}{\partial t}\right)\,\frac{\partial^{2}u^{\varepsilon}}{\partial t^{2}}
  u^\varepsilon \to 0\,\,\,\mbox{as}\,\,n\to\infty
\end{eqnarray}
Note that on $\Omega_{T}\setminus A_{\varepsilon}$, we have
\begin{eqnarray}\label{BVestimate.BVinitialdata.ean10}
 \left|sg^{\prime}_{n}\left(\frac{\partial u^\varepsilon}{\partial t}\right)\,\frac{\partial^{2}u^{\varepsilon}}{\partial t^{2}}
  u^\varepsilon\right|\leq sg^{\prime}_{n}\left(\frac{\partial u^\varepsilon}{\partial t}\right)\left|\frac{\partial^{2}u^{\varepsilon}}
  {\partial t^{2}}u^\varepsilon\right|  
\end{eqnarray}
Let $(x_{0}^{\varepsilon},t_{0}^{\varepsilon})\in A^{\varepsilon}$ and $R^{\prime}> 0$ be a real number. 
We denote the open ball with center at $(x_{0}^{\varepsilon},t_{0}^{\varepsilon})$ and having radius $R'$ by 
$B((x_{0}^{\varepsilon},t_{0}^{\varepsilon}),\,R^{\prime})$.\\

For each $n\in\N$, denote 
\begin{eqnarray}\label{mollifier.bv.eqn2}
 C_{n}^{\varepsilon}:=\left\{(x,t)\in\Omega_{T}\smallsetminus A^{\varepsilon}:\,\,0<\left|\frac{\partial u^{\varepsilon}}{\partial t}(x,t)
 \right|\leq\frac{1}{n}\right\}.
\end{eqnarray}
Observe that for each $n\in\mathbb{N}$,   we have $C_{n+1}^{\varepsilon}\subseteq C_{n}^{\varepsilon}$.
Since $\Omega_{T}$ is bounded, there exists $n_{0}\in\mathbb{N}$ be such that for all $n\geq n_{0}$, the following inclusion holds:
$$C_{n}^{\varepsilon}\subset B\left((x_{0}^{\varepsilon},t^{\varepsilon}\right),\frac{n}{2})$$ 
Define a function $\rho\in C^{\infty}_{0}(\mathbb{R}^{d+1})$   by
\begin{eqnarray}\label{mpllifier.bv.eqn11}
\rho(x,t) := \left\{\def\arraystretch{1.2}%
\begin{array}{@{}c@{\quad}l@{}}
  k_\varepsilon\,\exp\left(-\frac{1}{1-|(x,t)-(x_{0}^{\varepsilon},t_{0}^{\varepsilon})|^{2}}\right), & 
  \text{if $x\in B((x_{0}^{\varepsilon},t^{\varepsilon}_{0}),\,1)$}\\
  0,\,\, & \text{if\,\,$(x,t)\notin B((x_{0}^{\varepsilon},t_{0}^{\varepsilon}),\,1)$},
\end{array}\right.
\end{eqnarray}
where the constant $k_\varepsilon$ is chosen so that
\begin{eqnarray}\label{mollifier.bv.eqn4}
\int_{\mathbb{R}^{d+1}}\rho(x,t)\,\,dx=1.
\end{eqnarray}
Denote the sequence of mollifiers $\rho_{n}:\mathbb{R}^{d+1}\to\mathbb{R}$ by 
\begin{eqnarray}\label{mollifier.bv.eqn5}
\rho_{n}(x,t) := \left\{\def\arraystretch{1.2}%
\begin{array}{@{}c@{\quad}l@{}}
  k_\varepsilon\,n^{d+1}\,\exp\left(-\frac{n^{2}}{n^{2}-|(x,t)-(x_{0}^{\varepsilon},t_{0}^{\varepsilon})|^{2}}\right), & 
  \text{if $(x,t)\in B((x_{0}^{\varepsilon},t_{0}^{\varepsilon}),\,n)$}\\
  0,\,\, & \text{if\,\,$(x,t)\notin B((x_{0}^{\varepsilon},t_{0}^{\varepsilon}),\,n)$}
\end{array}\right.
\end{eqnarray}
Since for each $n\geq n_{0}$, we have 
\begin{eqnarray}\label{mollifier.bv.eqn6}
 sg_{n}^{\prime}\left(\frac{\partial u^{\varepsilon}}{\partial t}\right)= \left\{\def\arraystretch{1.2}%
  \begin{array}{@{}c@{\quad}l@{}}
   n & \text{if $0\leq\left|\frac{\partial u^{\varepsilon}}{\partial t}\right|\leq\frac{1}{n}$},\\
    0\,\, & \text{if\,\,$\left|\frac{\partial u^{\varepsilon}}{\partial t}\right|>\frac{1}{n},$}
    \end{array}\right.
\end{eqnarray}
therefore we compute
\begin{eqnarray}\label{mollifier.bv.eqn7}
\int_{\Omega_{T}\smallsetminus A^{\varepsilon}}\,sg_{n}^{\prime}\left(\frac{\partial u^{\varepsilon}}{\partial t}\right)\left|
\frac{\partial^{2}u^{\varepsilon}}{\partial t^{2}}u^\varepsilon\right|\,\,dx\,dt\hspace*{1.7in}\nonumber\\
=\int_{B((x_{0}^{\varepsilon},t_{0}^{\varepsilon}),\frac{n}{2})}\chi_{\Omega_{T}\smallsetminus A^{\varepsilon}}\,sg_{n}^{\prime}
\left(\frac{\partial u^{\varepsilon}}{\partial t}\right)\left|
\frac{\partial^{2}u^{\varepsilon}}{\partial t^{2}}u^\varepsilon\right|\,dx\hspace*{0.4in}\nonumber\\
= \int_{B((x_{0}^{\varepsilon},t_{0}^{\varepsilon}),\frac{n}{2})}\chi_{\Omega_T\smallsetminus A^{\varepsilon}}(x)\,\frac{sg_{n}^{\prime}
\left(\frac{\partial u^{\varepsilon}}{\partial t}\right)}{\rho_{n}((x,t))}\rho_{n}(x,t)\,\left|
\frac{\partial^{2}u^{\varepsilon}}{\partial t^{2}}u^\varepsilon\right|\,dx\,dt\nonumber\\
{}
\end{eqnarray}
For all $n\geq n_{0}$ and for $(x,t)\in B((x_{0}^{\varepsilon},t_{0}^{\varepsilon}),\,\frac{n}{2})$, we have
\begin{eqnarray}\label{mollifier.bv.eqn8}
 \frac{sg_{n}^{\prime}\left(\frac{\partial u^{\varepsilon}}{\partial t}\right)}{\rho_{n}(x,t)}:= \left\{\def\arraystretch{1.2}%
  \begin{array}{@{}c@{\quad}l@{}}
   \frac{1}{k_\varepsilon\,n^{d}}\,e^{\frac{n^{2}}{n^{2}-|(x,t)-(x_{0}^{\varepsilon},t_{0}^{\varepsilon})|^{2}}} & 
   \text{if $(x,t)\in B((x_{0}^{\varepsilon},t_{0}^{\varepsilon}),\,\frac{n}{2})\cap C_{n}^{\varepsilon}$},\\
    0\,\, & \text{if\,\,$(x,t)\notin B((x_{0}^{\varepsilon},t_{0}^{\varepsilon}),\,\frac{n}{2})\cap C_{n}^{\varepsilon}.$}
    \end{array}\right.
\end{eqnarray}
For $(x,t)\in  B((x_{0}^{\varepsilon},t_{0}^{\varepsilon}),\,\frac{n}{2})$, we obtain
$$\left| \frac{sg_{n}^{\prime}\left(\frac{\partial u^{\varepsilon}}{\partial t}\right)}{\rho_{n}(x,t)}\right|\leq\frac{1}{k_\varepsilon\,n^{d}}
\,e^{\frac{4}{3}}$$
\vspace{0.1cm}\\
Since $n\in\mathbb{N}$, the integrand on the last line of  \eqref{mollifier.bv.eqn7} is dominated by 
$$\frac{1}{k_\varepsilon}\,e^{\frac{4}{3}}\,\rho_{n}(x,t)\,\left|
\frac{\partial^{2}u^{\varepsilon}}{\partial t^{2}}u^\varepsilon\right|,$$
which is integrable on $\Omega_{T}$ as $u^{\varepsilon}\in C^{4+\beta,\frac{4+\beta}{2}}(\overline{\Omega_{T}})$ and 
$$\int_{\mathbb{R}^{d+1}}\rho_{n}(x,t)\,dx=1.$$
Therefore an application of dominated convergence theorem gives \eqref{BVestimate.BVinitialdata.ean8}. We now pass to the limit 
in the second term on RHS of  \eqref{BVestimate.BVinitialdata.ean7}. Note that 
\begin{eqnarray}\label{mollifier.bv.eqn8AA}
sg_{n}\left(\frac{\partial}{\partial t}u^\varepsilon(x,T)\right)u^{\varepsilon}(x,T)- sg_{n}\left(\frac{\partial}{\partial t}
u^\varepsilon(x,0)\right)u_{0\varepsilon}(x)\to sg\left(\frac{\partial}{\partial t}u^\varepsilon(x,T)\right)u^{\varepsilon}(x,T)
\nonumber\\- sg\left(\frac{\partial}{\partial t}
u^\varepsilon(x,0)\right)u_{0\varepsilon}(x)\,\,\mbox{as}\,\,n\to\infty\nonumber\\
{}
\end{eqnarray}
 and 
 $$\left|sg_{n}\left(\frac{\partial}{\partial t}u^\varepsilon(x,T)\right)u^{\varepsilon}(x,T)- sg_{n}\left(\frac{\partial}{\partial t}
u^\varepsilon(x,0)\right)u_{0\varepsilon}(x)\right|\leq 2\,\|u_{0}\|_{L^{\infty}(\Omega)},$$
which is integrable as $\mbox{Vol}(\Omega)<\infty$. Therefore
an application of dominated convergence theorem gives
\begin{eqnarray}\label{BVestimate.BVinitialdata.ean9}
 \int_{\Omega}\left(sg_{n}\left(\frac{\partial}{\partial t}u^\varepsilon(x,T)\right)u^{\varepsilon}(x,T)- sg_{n}\left(\frac{\partial}{\partial t}
u^\varepsilon(x,0)\right)u_{0\varepsilon}(x)\right)\,dx\,dt\nonumber\\ \to
 \int_{\Omega}\left(sg\left(\frac{\partial}{\partial t}u^\varepsilon(x,T)\right)u^{\varepsilon}(x,T)- sg\left(\frac{\partial}{\partial t}
u^\varepsilon(x,0)\right)u_{0\varepsilon}(x)\right)\,dx\,dt,\,\,\mbox{as}
 \,\,n\to\infty\nonumber\\
 {}
\end{eqnarray}
Using equations \eqref{BVestimate.BVinitialdata.ean8},\eqref{BVestimate.BVinitialdata.ean9} in \eqref{BVestimate.BVinitialdata.ean7}, we have
\begin{eqnarray}\label{BVestimate.BVinitialdata.ean10}
\int_{\Omega_{T}}\left|\frac{\partial u^\varepsilon}{\partial t}\right|\,dx\,dt=\int_{\Omega}
\left(sg\left(\frac{\partial}{\partial t}u^\varepsilon(x,T)\right)u^{\varepsilon}(x,T)- sg\left(\frac{\partial}{\partial t}
u^\varepsilon(x,0)\right)u_{0\varepsilon}(x)\right)\,dx\nonumber\\
{}
\end{eqnarray}
Applying Theorem \ref{Compactness.lemma.1} and using equations \eqref{BVestimate.BVinitialdata.ean1} and 
\eqref{BVestimate.BVinitialdata.ean10} in \eqref{BVestimate.BVinitialdata.ean5}, we get
\begin{eqnarray}\label{BVestimate.BVinitialdata.ean11}
 \varepsilon\int_{\Omega_{T}}\,B(u^{\varepsilon})\left|\Delta u^\varepsilon\right|\,dx\,dt &\leq& \|B^{\prime}\|_{L^{\infty}(I)}
 \frac{\mbox{Vol}(\Omega)}{2r}\|u_{0}\|_{L^{\infty}(\Omega)}^{2} + \displaystyle\max_{1\leq j\leq d}\left(\displaystyle\sup_{y\in I}
 \left|f_{j}^{\prime}(y)\right|\right)\,TV_{\Omega}(u_{0})\nonumber\\
 && + 2\,\|u_{0}\|_{L^{\infty}(\Omega)}\mbox{Vol}(\Omega).
\end{eqnarray}
Taking absolute value on both sides of \eqref{BVestimate.BVinitialdata.ean4} and integrating over $\Omega_{T}$, we have
\begin{eqnarray}\label{BVestimate.BVinitialdata.ean11}
 \int_{\Omega_{T}}\left|\frac{\partial u^{\varepsilon}}{\partial t}\right|\,dx\,dt &\leq& 2\|B^{\prime}\|_{L^{\infty}(I)}
 \frac{\mbox{Vol}(\Omega)}{2r}\|u_{0}\|_{L^{\infty}(\Omega)}^{2} + 2\displaystyle\max_{1\leq j\leq d}\left(\displaystyle\sup_{y\in I}
 \left|f_{j}^{\prime}(y)\right|\right)\,TV_{\Omega}(u_{0})\nonumber\\
 && + 2\,\|u_{0}\|_{L^{\infty}(\Omega)}\mbox{Vol}(\Omega).
\end{eqnarray}
\textbf{Step 3:} Denote the total variation of $u^{\varepsilon}$ by $TV_{\Omega_{T}}(u^{\varepsilon})$. Since for each $\varepsilon > 0$,
$u^{\varepsilon}\in H^{1}(\Omega_{T})$, the total variation of $u^{\varepsilon}$ is given by
\begin{eqnarray}\label{B.BVestimate28}
 TV_{\Omega_{T}}(u^{\varepsilon})= \left\|\frac{\partial u^{\varepsilon}}{\partial t}\right\|_{L^{1}(\Omega_{T})} +  
 \|\nabla u^{\varepsilon}\|_{\left(L^{1}(\Omega_{T})\right)^{d}}.
\end{eqnarray}
Using equations \eqref{BVestimate.BVinitialdata.ean1} and \eqref{BVestimate.BVinitialdata.ean2}, we get that 
$\left(TV_{\Omega_{T}}(u^{\varepsilon})\right)_{\varepsilon\geq 0}$ is bounded. Applying the fact that 
$BV(\Omega_{T})\cap L^{1}(\Omega_{T})$ is compactly imbedded in $L^{1}(\Omega_{T})$ \cite{MR1304494}, we get the existence of a  
subsequence $(u^{\varepsilon_{k}})$ and a function $u\in L^{1}(\Omega_{T})$ such that $u^{\varepsilon_{k}}\to u$ in $L^{1}(\Omega_{T})$ 
as $k\to\infty$. We still denote the subsequence by $(u^{\varepsilon})$. Since $u^{\varepsilon}\to u$ in $L^{1}(\Omega_{T})$ as  
$\varepsilon\to 0$, there exists a further subsequence $(u^{\varepsilon_{k}})$ of $\left(u^{\varepsilon}\right)$ such that we have 
\eqref{BVestimateeqn25} and \eqref{BVestimateeqn25.b}.
\vspace{0.3cm}\\
\subsection{BV estimate with initial data $u_{0}$ in $W^{1,1}_{0}(\Omega)\cap C(\overline{\Omega})$}
Following the proof of a result(Lemma 7.1) \cite[p.47]{Ramesh},  we conclude the following result. 
\begin{lemma}\label{paper2.initialdata.lemma1}
 Let $u_{0}\in W^{1,1}_{0}(\Omega)\cap C(\overline{\Omega})$. Then there exists a sequuence $\left(u_{0\varepsilon}\right)$ 
 in $\mathcal{D}(\Omega)$ such that the following properties hold.
 \begin{enumerate}
  \item As $\varepsilon\to 0$, we have
  \begin{eqnarray}\label{initialdata.h1.eqn1}
   u_{0\varepsilon}\to u_{0}\,\,\mbox{in}\,\,W^{1,1}(\Omega)
  \end{eqnarray}
\item For all $\varepsilon >0$, there exists a constant $A> 0$ such that 
\begin{eqnarray}\label{initialdata.h1.eqn2}
 \|u_{0\varepsilon}\|_{L^{\infty}(\Omega)}\leq A
\end{eqnarray}
 \end{enumerate}
\end{lemma}
Following exactly the same proof of Theorem \ref{paper3.BVestimate.Maintheorem}, using Hypothsis F and 
Theorem \ref{Compactness.lemma.1D}, we conclude that
\begin{theorem}\label{paper3.BVestimate.MaintheoremA}
 Let $f,\,B,\,u_{0}$ and $u_{0\varepsilon}$ satisfy Hypothesis F and $\left(u^{\varepsilon}\right)$ be as in 
 Lemma \ref{regularized.BVestimates.theorem3BB}.
 Then
 \begin{enumerate}
  \item for all $\varepsilon> 0$, there exists a constant $C> 0$ such that
  \begin{eqnarray}\label{BVestimate.BVinitialdata.ean1AB}
   \|\nabla u^{\varepsilon}\|_{\left(L^{1}(\Omega_{T})\right)}\leq C
  \end{eqnarray}
\item for all $\varepsilon > 0$, the following inequality  
\begin{eqnarray}\label{BVestimate.BVinitialdata.ean2AB}
\left\|\frac{\partial u^{\varepsilon}}{\partial t}\right\|_{L^{1}(\Omega_{T})}&\leq& 2\|B^{\prime}\|_{L^{\infty}(I)}
 \frac{\mbox{Vol}(\Omega)}{2r} A^{2} + 2\displaystyle\max_{1\leq j\leq d}\left(\displaystyle\sup_{y\in I}
 \left|f_{j}^{\prime}(y)\right|\right)\,C\nonumber\\
 && + 2\, A\mbox{Vol}(\Omega) 
\end{eqnarray}
holds.
 \end{enumerate}
 Further there exists a subsequence $(u^{\varepsilon_{k}})$ of $(u^{\varepsilon})$, and a function $u\in L^{1}(\Omega_{T})$ such that 
\begin{eqnarray}\label{BVestimateeqn25AB}
u^{\varepsilon_{k}}\to u\,\,\mbox{in}\,L^1(\Omega_T),\\
 u^{\varepsilon_{k}}\to u\,\,{\it a.e.}\,(x,t)\in \Omega_{T}
\end{eqnarray}
as $k\to\infty$.
\end{theorem}
\vspace{0.3cm}
\textbf{Proof of Theorem \ref{paper3.BVestimates.theorem2}:} The proof of Theorem \ref{paper3.BVestimates.theorem2} follows from the proof 
of Theorem 1.2 from \cite[p.40]{Ramesh}.\\
\vspace{0.2cm}
\textbf{Proof of Theorem \ref{paper3.BVestimates.theorem3}:} The proof of Theorem \ref{paper3.BVestimates.theorem3} follows from the proof 
of Theorem 1.3 from \cite[p.47]{Ramesh}.
\newpage
\bibliographystyle{amsplain}

\end{document}